\documentclass[11pt,leqno]{article}
\usepackage{amssymb,latexsym}

\def\del{\delta } \def\Del{\Delta } 
\def\gam{\gamma }

\def\sig{\sigma }

\def\lam{\lambda }

\def\part{\partial}

\def\var{\varphi}

\def\om{\omega }
\def\Om{\Omega}

\def\PP{\mathbf P}
 \newcommand{\begeq}{\begin{equation}}
\newcommand{\fineq}{\end{equation}}
 
 \newcommand{\ov}{\overline}
 \newcommand{\noi}{\noindent}

\font\piccola=cmr9 at 9pt
 \relax

\input prepictex
\input pictex
\input postpictex

\begin{document}

\begin{center}{\large\bf Error analysis for quadratic spline quasi-interpolants \\ on non-uniform
criss-cross triangulations \\ of bounded rectangular domains\\
Version 25/1/06} \end{center} \medskip

\begin{center} 
{\sc Catterina Dagnino} \\ Dipartimento di Matematica dell'Universit\`a di Torino \\ Via Carlo Alberto 10, 10123 Torino, Italy \\ {\it email: catterina.dagnino@unito.it}
 \end{center}

\begin{center}
{\sc Paul Sablonni\`ere} \\ INSA de Rennes,
20 avenue des Buttes de Co\"esmes, \\ CS 14315, 35043 Rennes Cedex, France. \\
 {\it email: Paul.Sablonniere@insa-rennes.fr} \end{center} \vspace{15pt}
  
\vspace{1cm}

\noindent{\bf Abstract.} 
{\it Given a non-uniform criss-cross partition of a rectangular domain $\Omega$, we analyse the error between a function $f$ defined on $\Omega$ and two types of
$C^1$-quadratic spline quasi-interpolants (QIs) obtained as linear combinations of
B-splines with discrete functionals as coefficients. The main novelties  are the
facts that supports of B-splines are contained in $\Omega$ and that data
sites also lie inside or on the boundary of $\Omega$.   Moreover, the infinity norms
of these QIs are small and do not depend on the triangulation: as the two QIs are
exact on quadratic polynomials, they give the optimal approximation order for smooth functions.
Our analysis is done for $f$ and its partial derivatives of the first and
second orders and a particular effort has been made in order to give the best possible
error bounds in terms of the smoothness of $f$ and of the mesh ratios of the
triangulation.}

\vspace {.5cm} 

\vspace {.3cm} 

\vspace {.3cm} 

\section{Introduction} 

Given a non-uniform criss-cross partition of a rectangular domain $\Omega$, we
analyse the error between a function $f$ defined on $\Omega$ and two
$C^1$ {\sl quadratic spline quasi-interpolants} (abbr. QIs), denoted $S_2$ and
$W_2^*$, obtained as linear combinations of B-splines with discrete coefficient functionals.
The first operator $S_2$ was described by the second author in \cite{s1}\cite{s4} and the second
one $W_2^*$ is a slight modification of the operator $W_2$ introduced by
Chui and Wang in \cite{cw1}, and also studied by Chui and He in \cite{ch}, Wang and Lu in \cite{wl} and  by the first author in  \cite{dl1}, \cite{dl2}.
\\
With respect to previous papers, we note the following facts :  we introduce  B-splines  with {\sl supports  contained  in} $\Omega$ and {\sl data sites}
lying {\sl inside or on the boundary of} $\Omega$, so we do not need extra values outside the domain.  
This can be useful in certain practical problems where these data are not available. Moreover, we show that the infinity norms of
these QIs are small and {\sl do not depend on the triangulation}. As they are exact on the space $\PP_2$ of quadratic polynomials, it is well known that they give the {\sl optimal approximation order} for  smooth functions.\\
Another important and very useful property of QIs is that the construction
of these operators {\sl do not need the solution of any system of equations}. It is
particularly attractive in the bivariate case where the number of data sites can be huge
in practice.
\\
Though the QIs do not interpolate $f$ at data sites, it can be observed that errors are quite small at that points. Actually,  a superconvergence phenomenon can often be observed at some specific points. Moreover, the global
behaviours of QIs and of their derivatives is quite close to those of the function $f$ (see e.g.\cite{fs}).
\\
Our error analysis is done for $f$ and its partial derivatives of the first and second orders and a particular effort has been made in order to obtain sharp error bounds in terms of
the smoothness of $f$ and of the characteristics of the triangulation, in particular local mesh ratios. Such a program can be
developed thanks to the good properties of quadratic B-splines described in \cite{s3}. It is true that we do not get the {\sl best} error constants, which is a rather technical task, however, we obtain a reasonable  order of magnitude of these constants. This can be useful in the practical applications that we want to develop elsewhere.
\\
Here is an outline of the paper: in Section 2, we recall the main definitions on the
B-splines on criss-cross triangulations that we use in the definition of quasi-interpolants . In Section 3,  we describe the two quadratic spline QIs.
In Section 4, we give error estimates of the infinity norms of $f-Q$, where $Q=S_2$ or
$W_2^*$, when $f\in C^s(\Omega)$, with $0\le s\le 3$.
In Section 5, we give error estimates on  first derivatives $\Vert D^{r,s}(f-Q)\Vert_{\infty},\, r+s=1$, in $\Omega$, and on second derivatives $\Vert D^{r,s}(f-Q)\Vert_{\infty},\, r+s=2$, inside triangular cells of the triangulation, since $Q$ is only $C^1$. They are expressed in terms of moduli of smoothness with respect to the length $h/2$, where $h$ is the maximal steplength of the given partition of the domain.


\section{Quadratic B-splines  on a bounded rectangle}

In this Section, we first introduce $C^1$ quadratic B-splines generating the spline space in which we approximate functions. Then, in the following section, we will define the two quasi-interpolants $S_2$ and $W_2^*$.

Let $\Om = [a,b] \times [c,d]$ be a rectangle decomposed into $mn$ subrectangles by the two partitions
$$
X_m=\{ x_i,0\le i\le m\}, \quad Y_n=\{ y_j, 0\le j \le n\},
$$
respectively of the segments $I=[a,b]=[x_0,x_m]$ and $J=[c,d] = [y_0,y_n]$. We also introduce the double knots $x_{-1}=x_0, \,y_{-1}=y_0, \,x_{m+1}=x_m, \,y_{n+1}=y_n$.

The so-called {\sl criss-cross triangulation} ${\cal T}_{mn}$ of $\Omega$ is defined by drawing the two diagonals in each subrectangle $\Om_{ij}=[x_{i-1},x_i] \times [y_{j-1},y_j]$. 
We need the two following sets of indices:

$$
\begin{array}{ll}
{\cal K}_{mn}=\{ (i,j) : 0\le i \le m+1, \quad 0\le j \le n+1\}, \\
\\
\widehat {\cal K}_{mn} = \{(i,j) : 1\le i \le m, \quad \> 1\le j \le n\}. 
\end{array}
$$

We set $h_i=x_i-x_{i-1},\; k_j=y_j-y_{j-1},\; s_i=\frac 12 (x_{i-1}+x_i),\;t_j=\frac  12 (y_{j-1}+y_j)$
for $(i,j)\in {\cal K}_{mn}$. 
We denote by $\{A_{i,j}=(x_i,y_j),\,0\le i\le m, \,0\le j\le n\}$ the set of vertices of subrectangles and by $\{M_{i,j}=(s_i,t_j),\,(i,j)\in {\cal K}_{mn} \}$ the set of their centers, of  midpoints of boundary subintervals and of vertices of $\Omega$.

Let ${\cal B}_{mn}:=\{ B_{ij} , (i,j) \in {\cal K}_{mn}\}$  be the collection of $(m+2)(n+2)$ B-splines generating the space ${\cal S}_2({\cal T}_{mn})$ of all $C^1$ piecewise quadratic functions on the criss-cross triangulation ${\cal T}_{mn}$, associated with the partition $X_m\times Y_n$ of the domain $\Omega$. There are $mn$ B-splines associated with the set of indices $\widehat {\cal K}_{mn}$, whose restrictions to the boundary $\Gamma$ of $\Omega$ are equal to zero. They were also introduced in \cite{csw}\cite{cw1}\cite{cw2}. To the latter, we add $2m+2n+4$ {\em boundary B-splines} whose restrictions to $\Gamma$ are univariate quadratic B-splines. Their set of indices is
$$
\widetilde {\cal K}_{mn} : = \{(i,0),(i,n+1),0\le i \le m+1; (0,j),(m+1,j), 0\le j \le n+1\}.
$$ 

The BB (=Bernstein-B\'ezier)-coefficients of inner B-splines  $\{B_{ij},\; 2\le i \le m-1$, $2\le j \le n-1 \}$ are given in \cite{s1}. The other ones can be found in the technical report \cite{s3} and in \cite{s5}. The B-splines  are positive and form a partition of unity (blending system). The boundary B-splines are {\em linearly independent} as the univariate ones. But the inner B-splines are {\em linearly dependent}, the dependence relationship being:
$$
\sum_{(i,j) \in \hat {\cal K}_{mn}} (-1)^{i+j}h_ik_j B_{ij}=0.
$$

Although ${\cal B}_{mn}$ is not a basis of ${\cal S}_2({\cal T}_{mn})$, this fact has no influence on the definition and properties of QIs. The support of $B_{ij}$ is denoted by $\Sigma_{ij}$: for inner B-splines, it is a non-uniform octagon.
The set ${\cal B}_{mn}$ can also be defined in the following way.  Define the extended partitions 
$$
\ov X_m= X_m  \cup \{ \ov x_{-2}, \ov x_{-1}, \ov x_{m+1}, \ov x_{m+2}\} 
$$
and
$$
 \ov Y_n=Y_n  \cup \{\ov y_{-2}, \ov y_{-1}, \ov y_{n+1}, \ov y_{n+2}\},
$$
 where $\ov x_{-2} < \ov x_{-1} < x_{0},\quad x_{m} < \ov x_{m+1} < \ov x_{m+2},\quad \ov y_{-2} < \ov y_{-1} < y_{0},\quad y_{n} < \ov y_{n+1} < \ov y_{n+2}$,  and the corresponding criss-cross triangulation $\ov{\cal T}_{mn}$. We also 
 put $\bar h_0=x_0-\ov x_{-1},\, \bar h_{m+1}=\ov x_{m+1}-x_m,\, \bar k_0=y_0-\ov y_{-1},\, \ov k_{n+1}=\ov y_{n+1}-y_n.$

We consider the collection $\ov{\cal B}_{mn}:=\{\ov{B}_{ij} , (i,j) \in {\cal K}_{mn}\}$ of the  $(m+2)(n+2)$ "classical" B-splines with octagonal supports  $\bar\Sigma_{ij}$ such that $\bar\Sigma_{ij}\cap int(\Omega)\not=\emptyset$ \cite{cw1}.

We note that $B_{ij}=\ov{B}_{ij}$ for inner B-splines. Using the BB-coefficients of both families ${\cal B}_{mn}$ and $\ov{\cal B}_{mn}$, one can derive the expressions of the new boundary B-splines in function of "classical" B-splines. 
For this purpose, we need the following notations, for $2\le i\le m$ and $2\le j\le n$ :
$$
\sig_i=\frac  {h_i} {h_{i-1}+h_i},  \;\;\sig_i' = \frac  {h_{i-1}} {h_{i-1}+h_i} = 1-\sig_i \> , \quad
\tau_j=\frac  {k_j} {k_{j-1}+k_j}, \;\; \tau_j' = \frac  {k_{j-1}} {k_{j-1}+k_j} = 1-\tau_j \> . 
$$
In addition, we need  the particular values :
$$
\ov  \sig_1=\frac  {h_1} {\ov  h_{0}+h_1},  \;\; \ov \sig'_{m+1}= \frac  {h_{m}} {h_{m}+ \ov  h_{m+1}}, \quad
\ov \tau_1=\frac  {k_1} {\ov k_{0}+k_1}, \;\; \ov \tau'_{n+1}= \frac  {k_{n}} {k_{n}+ \ov  k_{n+1}},
$$
and  $\sigma_1=\sigma'_{m+1}=\tau_1=\tau'_{n+1}=1$, whence $\sigma'_1=\sigma_{m+1}=\tau'_1=\tau_{n+1}=0$, since $h_0=h_{m+1}=k_0=k_{n+1}=0$.\\

The {\sl first boundary layer} of B-splines along the horizontal edge $A_{00}A_{m0}$ is defined by
$$
\begin{array}{ll}
B_{00}={\frac{1}{\ov  \sig_1 \ov  \tau_1}}\ov{B}_{00},\quad B_{10}={\frac{1}{\ov  \tau_1}}(\ov{B}_{10}- {\frac{\ov  \sig_1'}{\ov  \sig_1}}\ov{B}_{00}),\quad
B_{i0}={\frac{1}{\ov  \tau_1}}\ov{B}_{i0}, \quad 2\le i \le m-1,\\
\\
 B_{m0}={\frac{1}{\ov  \tau_1}}(\ov{B}_{m,0}- {\frac{\ov  \sig_{m+1}}{\ov  \sig'_{m+1}}}\ov{B}_{m+1,0}), \quad B_{m+1,0}={\frac{1}{\ov  \sig'_{m+1} \ov  \tau_1}}\ov{B}_{m+1,0}.
\end{array}
$$
\\
In the same way we obtain, along the vertical edge $A_{00}A_{0n}$,
$$
\begin{array}{ll}
B_{01}={\frac{1}{\ov  \sig_1}}(\ov{B}_{01}- {\frac{\ov  \tau_1'}{\ov  \tau_{1}}}\ov{B}_{0,0}), \quad
B_{0j}={\frac{1}{\ov  \sig_1}}\ov{B}_{0j}, \quad 2\le j \le n-1,\\
\\
 B_{0n}={\frac{1}{\ov  \sig_1}}(\ov{B}_{0n}- {\frac{\ov  \tau_{n+1}}{\ov  \tau'_{n+1}}}\ov{B}_{0,n+1}), \quad
B_{0,n+1}=\frac{1}{\ov  \sig_1 \ov  \tau'_{n+1}}\ov{B}_{0,n+1}.
\end{array}
$$
\\
Similar formulas hold for boundary B-splines along  the edges $A_{m0}A_{mn}$ and $A_{0n}A_{mn}$ : \\
$$
\{B_{i,n+1},0\le i\le m+1\} \;\;{\rm and}\;\; \{B_{m+1,j},0\le j\le n+1\}.
$$

The restrictions of all these B-splines to the boundary of $\Om$ are classical univariate quadratic B-splines.
The {\sl second boundary layer} of B-splines along the horizontal edge $A_{00}A_{m0}$ is defined by
$$
\begin{array}{ll}
B_{11}=\ov{B}_{11}-{\frac{\ov \tau_1'}{\ov \tau_{1}}}\ov{B}_{1,0}-{\frac{\ov \sig_1'}{\ov\sig_{1}}}\ov{B}_{01}+{\frac{\ov \sig'_1\ov\tau'_1}{\ov\sig_1\ov \tau_1}}\ov{B}_{00},
\end{array}
$$
and similar formulas for $B_{m1}, B_{1n}$ and  $B_{mn}$.\\

Finally we define the second layer of B-splines along the vertical edge $A_{00}A_{0n}$,
$$
B_{i1}=\ov{B}_{i1}-{\frac{\ov  \tau'_1}{\ov  \tau_{1}}}\ov{B}_{i,0},\quad 2\le i\le m-1,\quad
B_{1j}=\ov{B}_{1j}-{\frac{\ov  \sig'_1}{\ov  \sig_{1}}}\ov{B}_{0,j}, \quad 2\le j\le n-1,
$$
and similar formulas for the collections:
$$
\{B_{in}, \; 2\le i\le m-1\} \quad {\rm and} \quad \{B_{mj},\; 2\le j\le n-1\},
$$
with the ratios $\,\displaystyle\frac{\ov  \tau_{n+1}}{\ov  \tau_{n+1}'}\,$ and $\,\displaystyle\frac{\ov  \sig_{m+1}}{\ov  \sig_{m+1}'}\,$ instead of $\,\displaystyle\frac{\ov  \tau'_1}{\ov \tau_1}\,$and $\,\displaystyle\frac{\ov  \sig'_1}{\ov  \sig_1}\,$ respectively, in the formula defining $B_{11}.$\\

{\bf Remark : } note that many coefficients can be simplified, for example
$$
\frac{\bar \sig'_1}{\bar \sig_{1}}=\frac{\ov h_0}{h_1},\quad \frac{\bar \sig_{m+1}}{\bar \sig'_{m+1}}=\frac{\ov h_{m+1}}{h_m},\quad \frac{\bar \tau'_1}{\bar\tau_1}=\frac{\ov k_0}{k_1}, \quad \frac{\bar \tau_{n+1}}{\bar \tau'_{n+1}}=\frac{\ov k_{n+1}}{k_n}.
$$

Error analyses given in sections 3 and 4 below are based on the Bernstein B\'ezier representation of B-splines on the triangulation. The associated techniques are described e.g. in \cite{c}, \cite{f1}, \cite{p1}.


\section{ Quasi-Interpolants exact on $\PP_2$}

We now define the two quadratic spline quasi-interpolants  $S_2$ and $W_2^*$ that we want to study. Moreover, we give uniform bounds on their infinity norms.


\subsection{The quasi-interpolant $S_2$}

For the definition of $S_2$, we need the notations $\sig_i$ and $\tau_j$ given above in Section 2. Then we define:
$$
a_i=-\frac  {\sig^2_i \sig_{i+1}'} {\sig_i + \sig_{i+1}'} , \;\;
c_i = - \frac  {\sig_i(\sig_{i+1}')^2} {\sig_i + \sig_{i+1}'},  \;\;
\ov a_j = \frac  {\tau^2_j \tau_{j+1}'} {\tau_j + \tau_{i+1}'}, \;\;
 \ov  c_j = - \frac {\tau_j(\tau_{j+1}')^2}{\tau_j + \tau_{j+1}'}, 
 $$
 $$
b_{ij}=1-(a_i+c_i+\ov a_j +\ov c_j ),
$$
with $a_0=c_0=a_{m+1}=c_{m+1}=\ov a_0=\ov c_0=\ov a_{n+1}=\ov c_{n+1}=0$ and $b_0=\ov b_0= b_{m+1}=\ov b_{n+1}=1$.\\

The data sites for $S_2$ are the $(m+2)(n+2)$ points of the set
 $$
 {\cal D}_{mn} : = \{ M_{i,j} = (s_i,t_j), (i,j) \in {\cal K}_{mn}\} \>.
 $$

The  quadratic spline quasi-interpolants $S_2$ \cite{s2}\cite{s4} is defined as follows:
$$
S_2f = \sum^{m+1}_{i=0} \sum^{n+1}_{j=0} \mu_{ij}(f) B_{ij}, 
$$
with coefficient functionals given by
\begin{equation}
 \mu_{ij}(f)= b_{ij}f (M_{i,j}) + a_if (M_{i-1,j}) + c_i f(M_{i+1,j})
 + \ov a_j f(M_{i,j-1}) + \ov c_j f(M_{i,j+1}). 
\end {equation}

It is exact on $\PP_2$ and its infinity norm is uniformly bounded {\sl independently of the triangulation} ${\cal T}_{mn}$ of the domain. Indeed, since
\begeq
|a_i|,|c_i|,|\ov a_j|,|\ov c_j|\le 1/2\quad {\rm and} \quad  |b_{ij}|\le 3,\label{due}
\fineq 
then it is clear that
$$
||S_2||_\infty \le 5. 
$$ 
\\
We notice that the number of data sites requested by $S_2$ is equal to
\begeq
N_S=mn+2m+2n+4.\label{quattro0}
\fineq 

\subsection{The quasi-interpolant $W_2^*$}

The second quasi-interpolant $W_2^*$ here analysed is a modification of the QI $W_2$ derived by Chui-Wang  \cite{cw1}, which is also exact for $\PP_2$.  The latter is defined in terms of classical B-splines  $\{\ov{B}_{ij}\}$ on the triangulation
 $\ov {\cal T}_{mn}$.
Given the values of a function $f$ at the $(m+3)(n+3)$ points $A_{ij}=(x_i,y_j) \> , \>  -1\le i \le m+1$, $-1 \le j \le n+1$ (among which those having one extra abscissa or ordinate are outside $\Omega$) and the $(m+2)(n+2)$ points $\ov M_{ij}$, intersections of the diagonals in the subrectangles with vertices $ A_{i-1,j-1},  A_{i,j}, A_{i-1,j}, A_{i,j-1}$ (among which a number also lay outside the domain), the Chui-Wang QI is defined by:
$$
W_2 f = \sum^{m+1}_{i=0} \sum^{n+1}_{j=0} \ov \mu_{ij}(f) \ov B_{ij},
$$
with coefficient functionals defined by
$$
\ov \mu_{ij} (f) = 2f(\ov M_{i,j}) - \frac  14 [f(A_{i-1,j-1}) + f(A_{i-1,j}) + f(A_{i,j-1}) + f(A_{i,j})].
$$

In that case, the number of data sites  is equal to
$$
N_W=2mn+3m+3n+9.
$$
Now if we set $x_{-2} = x_{-1} =x_0, \; x_{m+2} = x_{m+1} = x_m,\;y_{-2} =y_{-1} = y_0,\; y_{n+2}=y_{n+1}=y_n$ and if we use the  B-splines $B_{ij}$ defined in Section 2, with supports $\Sigma_{ij}$ included in the domain $\Om$, we can define the {\em modified} Chui-Wang QI as follows:
$$
W_2^* f = \sum^{m+1}_{i=0} \sum^{n+1}_{j=0} \mu_{ij}^*(f) B_{ij},
$$
where the coefficient functionals are: 
\begeq
\mu^*_{ij}(f) = 2f(M_{i,j}^*) -\frac  14 [ f(A_{i-1,j-1}^*) + f(A_{i-1,j}^*) + f (A_{i,j-1}^*) + f(A_{i,j}^*)] \>, \label{quattro}
\fineq 

with the new data points :
\begin{eqnarray}
A^*_{i,j}&=&A_{ij}, \quad {\rm for } \quad 0\le i\le m,\quad 0\le j\le n,\nonumber\\
A^*_{i,-1} &=&A_{i,0},\quad  A_{i,n+1}^* = A_{i,n},\quad -1\le i \le m+1,\label{cinque}\\ 
A^*_{-1,j} &=& A_{0,j},\quad A^*_{m+1,j} = A_{m,j}, \quad -1\le i \le m+1,\nonumber\\
M_{i,j}^* &=&  M_{i,j} \quad {\rm for } \quad 0\le i\le m+1,\quad 0\le j\le n+1.\nonumber
\end{eqnarray} 

The number of data sites requested by $W_2^*$ is equal to
\begeq
N_W^*=2mn+m+n+1,\label{sei0}
\fineq 
and they all lie inside  the  domain $\Om$ or on its boundary. \medskip 

From $|\mu_{ij}^*(f)| \le 3 ||f||_\infty$, we can immediately deduce:
$$
||W_2^*||_{\infty} \le 3 
$$
for all non-uniform triangulations ${\cal T}_{mn}$ of the domain $\Omega$. \\

We remark that both $S_2$ and $W_2^*$ are local schemes, because for $(x,y)\in\Omega$, the values $S_2f(x,y)$ and $W_2^*(x,y)$ only depend on those of $f$ in a neighbourhood of $(x,y) $. 

\begin{figure}[h!]
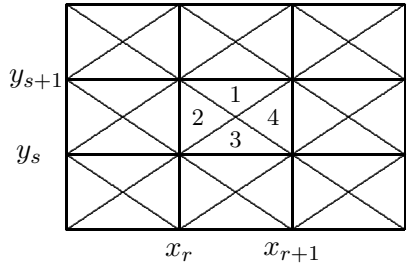

\caption{Four different kinds of cells in ${\cal T}_{mn}$.}
\label{fig4}
\end{figure}

If $0\le r\le m-1$ and $0\le s\le n-1$ are integers such that $x\in[x_{r}, x_{r+1}],\, y\in [ y_{s}, y_{s+1}]$, then $(x,y)$ will belong to one of the four triangular cells $T_\ell$ of ${\cal T}_{mn},\ \ell=1,2,3,4,$ labelled as in Fig. \ref{fig4}.

Each triangle $T_{\ell},\ \ell=1,2,3,4,$ is covered by exactly seven supports of B-splines $\Sigma_{ij}$. In Table 1 below,  we report the set $K(T_{\ell})$ of indices of such B-splines, as functions of $r$ and $s$, i.e. $K(T_{\ell})=\{(i,j)|\,\Sigma_{ij}\cap int(T_\ell)\not=\emptyset\}$.\\

Therefore, if $(x,y)\in T_{\ell},$ then:
 
$$
S_2f(x,y)=\sum_{(i,j)\in K(T_{\ell})}\mu_{ij}(f) B_{ij}(x,y)
$$
$$
W^*_2f(x,y)=\sum_{(i,j)\in K(T_{\ell})}\mu^*_{ij}(f) B_{ij}(x,y).
$$

\vskip5pt
$$
\beginpicture
\linethickness=.5pt
\setcoordinatesystem units <1cm,1cm>
\setplotarea x from 0 to 4.5, y from 0 to 3
\grid 3 3
\setlinear
\plot 0    0  4.5  3 /  
\plot 1.5  0  4.5  2 /  
\plot 3    0  4.5  1 /
\plot 0    1  3    3 /
\plot 0    2  1.5  3 /
\plot 0    3  4.5  0 /
\plot 1.5  3  4.5  1 /
\plot 3    3  4.5  2 /
\plot 0    2  3    0 /
\plot 0    1  1.5  0 /
\put {$x_r$} at 1.5  -.3
\put {$y_s$} at -.5    1
\put {$x_{r+1}$} at 3  -.3
\put {$y_{s+1}$} at -.4    2
\put {\piccola 1 } at  2.3  1.8
\put {\piccola 2 } at  1.8  1.5
\put {\piccola 3 } at  2.3  1.2
\put {\piccola 4 } at  2.8  1.5
\endpicture
$$
\vskip-15pt
\vskip5pt

\begin{center}
\begin{tabular}{|c|c|c|c|c|}
\hline
 & $T_1$ & $T_2$ & $T_3$ & $T_4$\\
\hline
\hline
& $r,\ s-1$ & $r-1,\ s-1$ & $r-1,\ s-1$ & $r,\ s-1$ \\
& $r-1,\ s$ & $r,\ s-1$ & $r,\ s-1$ & $r+1,\ s-1$ \\
& $r,\ s$ & $r-1,\ s$ & $r+1,\ s-1$ & $r-1,\ s$ \\
$i,j$ & $r+1,\ s$ & $r,\ s$ & $r-1,\ s$ & $r,\ s$ \\
& $r-1,\ s+1$ & $r+1,\ s$ & $r,\ s$ & $r+1,\ s$ \\
& $r,\ s+1$ & $r-1,\ s+1$ & $r+1,\ s$ & $r,\ s+1$ \\
& $r+1,\ s+1$ & $r,\ s+1$ & $r,\ s+1$ & $r+1,\ s+1$ \\
\hline
\hline
\end{tabular}
\end{center}
\vskip5pt
\centerline{\small Table 1.}


\section{Error analysis for functions}
 
In this section we analyse the errors $f - S_2f$ \ and \ $f-W_2^*f$ for $f\in C^s(\Omega)$, $0\le s \le 3.  $ \ We need to introduce the following notations:
\begin{eqnarray}
&& h = \max\{h_i,\,k_j\},\quad \delta= \min \{h_i,\,k_j\} ;\quad
||\cdot ||_{\infty,\Omega}= ||\cdot ||_\Omega= {\rm supremum \ norm \ over \ } \Omega; \nonumber\\
&&D^{\alpha}=D^{(\alpha_1,\alpha_2)} = \frac{\part^{\vert\alpha\vert}} {\part x^{\alpha_1} \part y^{\alpha_2}}\;\;{\rm with}\;\; \vert\alpha\vert=\alpha_1+\alpha_2\nonumber;\quad
 \om (D^sf, t)= \max\{\om (D^{\alpha}f, t), \vert\alpha\vert=s\};\nonumber\\ 
&& ||(x,y)||=(x^2+y^2)^{1/2}; \quad e_\alpha(x,y)=x^{\alpha_1}y^{\alpha_2}={\rm monomial\; of\; total \;degree} \;\vert\alpha\vert, \nonumber
\end{eqnarray}
where  the modulus of continuity of $\psi\in C(\Omega)$ is given by:
$$
\om(\psi, t)=\max \{ |\psi(M) -\psi(P)| ; \,M,P \in \Omega , ||MP|| \le t \}.
$$
We denote by $Q$ the generic quasi-interpolant defined by:
\begeq
Q f = \sum^{m+1}_{i=0} \sum^{n+1}_{j=0} \lam_{ij}(f) B_{ij}, \label{dieci} \fineq
where the coefficient  functionals are defined by $\lam_{ij}=\mu_{ij}$ (1) when $Q=S_2$, and  $\lam_{ij}=\mu_{ij}^*$ (4) when $Q=W_2^*$.\\


{\bf Theorem 1.} (Error bounds for continuous functions).
{\it There exists a constant $C_0>0$, with $C_0\le 20.5$ for $Q=S_2$ and $C_0 \le 12$ for $Q=W_2^*$, such that, for $f\in C(\Omega)$ }
$$
||f-Qf||_\Omega \le C_0\, \om (f, \frac12 h ).
$$
{\it Proof.} We consider some closed triangular cell $T$ of ${\cal T}_{mn}$, for which
$$
||f-Qf||_\Om = ||f-Qf||_{T}.
$$
$T$ is one of the four triangles depicted in Fig. \ref{fig4}. For the sake of simplicity we assume that         $T=T_3$.\
Since $Q$ reproduces $\PP_2$, for any $P \in T$, we can write $f(P)=\sum_{(i,j)\in K(T)} f(P) B_{ij}$.

For $Q=S_2$, we can write the following inequality:
\begin{eqnarray}
\mid (S_2f-f)(P)\mid&\le&\sum_{(i,j)\in K(T_3)}B_{ij}\,\{\vert b_{ij}\vert\vert f(M_{ij})-f(P)\vert+\vert a_i\vert\vert f(M_{i-1,j})-f(P)\vert+ \nonumber\\
&& \vert c_i\vert\vert f(M_{i+1,j})-f(P)\vert+\vert\overline a_j\vert\vert f(M_{i,j-1})-f(P)\vert+\vert\overline c_j\vert\vert f(M_{i,j+1})-f(P)\vert\} \nonumber
\end{eqnarray}\
Assuming that the origin lies at the midpoint of the lower edge of $T_3$, then this triangle can be decomposed into two equal subtriangles by the $y$-axis. By the symmetry of the problem, it is sufficient to consider the case when the point $P=(x,y)$ lies in the right triangle. Therefore the coordinates satisfy $0\le x+y \le \frac{h}{2}$. We shall now use the following simplified notations: there are seven B-splines whose supports intersect $int(T)$ and we denote their centres by $\{M_k,\,1\le k\le 7\}$, with $M_1=M_{r,s+1}$, $M_2=M_{r-1,s}$, $M_3=M_{r,s}$, $M_4=M_{r+1,s}$, $M_5=M_{r-1,s-1}$, $M_6=M_{r,s-1}$, $M_7=M_{r+1,s-1}$. Each central point $M_k$ has four neighbours $N_k,S_k,E_k,W_k$ (for North, South, East and West positions) involved in the coefficient functional $\mu_k$. The biggest constants being obtained for $k=1,2,5$, we only detail one of these cases, for example $k=1$.
Then, we  obtain the following upper bounds for the various distances involved in the majoration :
$$
\Vert PM_1\Vert\le\sqrt{10}\,\frac{h}{2}\le 4\,\frac{h}{2}, \quad\Vert PN_1\Vert\le  \sqrt{26}\,\frac{h}{2}\le 6\,\frac{h}{2},\quad\Vert PS_1\Vert\le \sqrt{2}\,\frac{h}{2}\le 2\,\frac{h}{2},$$
$$
\Vert PE_1\Vert\le \sqrt{13}\,\frac{h}{2}\le 4\,\frac{h}{2}, \;\;\Vert PW_1\Vert\le 3\sqrt{2}\,\frac{h}{2}\le 5\,\frac{h}{2}.
$$
Using inequalities (\ref{due}), we see that the coefficient of the B-spline $B_1=B_{r,s+1}$, whose support is centered at $M_1$ is first bounded above by
$$
3\omega( f,\Vert PM_1\Vert)+\frac12\left( \omega( f,\Vert PN_1\Vert)+\omega( f,\Vert PS_1\Vert)+\omega( f,\Vert PE_1\Vert)+\omega( f,\Vert PW_1\Vert) \right),
$$
then, using the above upper bounds on distances, we see that it is bounded above by :
$$
 [12+\frac12(6+2+4+5)] \,\omega(f, \frac{h}{2})=20.5\, \omega(f, \frac{h}{2}).
$$
Finally, since $\sum_{(i,j)\in K(T_3)}B_{ij}=1$, we obtain, for all $P\in T$ :
$$
\mid (S_2f-f)(P)\mid\le 20.5\, \omega(f, \frac{h}{2}).
$$
which proves that $\Vert f-S_2f \Vert \le 20.5\, \omega(f, \frac{h}{2})$.\\

Similarly, for $Q=W^*_2$, we can write the following inequality:
\begin{eqnarray}
\mid (W^*_2f-f)(P)\mid&\leq& \sum_{(i,j)\in K(T_3)}B_{ij}\,\{2\vert f(M_{i,j})-f(P)\vert+\frac14
\{ \vert f(A_{i,j})-f(P)\vert + \vert f(A_{i+1,j})-f(P)\vert\nonumber\\
&&+ \vert f(A_{i,j+1})-f(P)\vert + \vert f(A_{i+1,j+1})-f(P)\vert )\}.\nonumber
\end{eqnarray}\
We now compute upper bounds for the distances involved in the case when $(i,j)=(r,s+1)$. The central point $M_1=M_{r,s+1}$ has four neighbours  $NW_1, NE_1,SW_1, SW_1$ (for North-West, North-East, South-West and South-East positions) involved in the coefficient functional $\mu^*_1$, thus we obtain :
$$
\Vert PM_1\Vert\le\sqrt{10}\,\frac{h}{2}\le 4\,\frac{h}{2}, \quad 
\Vert PNW_1\Vert\le 2\sqrt{5}\, \frac{h}{2}\le 5\, \frac{h}{2} ,\quad 
\Vert PNE_1\Vert\le  \sqrt{17}\, \frac{h}{2}\le 5\, \frac{h}{2}, 
$$
$$
\quad \Vert PSW_1\Vert\le 2\sqrt{2}\, \frac{h}{2}\le 3\,\frac{h}{2},  \quad 
\Vert PSW_1\Vert\le \sqrt{5}\, \frac{h}{2}\le 3\,\frac{h}{2}.
$$
We see that the absolute value of the coefficient of the B-spline $B_1=B_{r,s+1}$ whose support is centered at $M_1$ is first bounded above by
$$
2\omega( f,\Vert PM_1\Vert)+\frac14\left( \omega( f,\Vert PNW_1\Vert)+\omega( f,\Vert PNE_1\Vert)+\omega( f,\Vert PSW_1\Vert)+\omega( f,\Vert PSE_1\Vert)\right),
$$
then, using the above upper bounds on distances, we obtain as upper bound :
$$
 [8+\frac14(5+5+3+3)] \,\omega(f, \frac{h}{2})=12 \, \omega(f, \frac{h}{2}),
$$
and, finally we obtain, for all $P\in T$ :
$$
\mid (W_2^*f-f)(P)\mid\le 12\, \omega(f, \frac{h}{2})
$$
which proves that $\Vert f-W_2^*f \Vert \le 12\, \omega(f, \frac{h}{2})$.
$\blacksquare$\\


\noi{\bf Theorem 2.} (Error bounds for $C^1$-functions). {\it  There exists
a constant $C_1>0$, with $C_1 \le 3$ for $Q=S_2$ and $C_1 \le 2$ for $Q=W_2^*$, such that, for $f\in C^1(\Om)$  :}
$$
||f-Qf||_\Om \le C_1 h\,\om (Df , h/2 ).
$$
\noi{\it Proof.} 
Let $q^*$ be the best approximation polynomial of $f$ in $\PP_1$ on the domain $\Omega$. Consider some closed triangular cell $T$ of $\mathcal{T}_{mn}$ in which we have
$$
\Vert f-q^*\Vert_T=\Vert f-q^*\Vert_\Omega.
$$
Take a point $(\xi,\eta)$ at the midpoint of the external edge of $T$ and let $q_1\in \PP_1$ be the linear Taylor polynomial of $f$ at that point :
\begeq
q_1(x,y)=f(\xi,\eta)+D^{(1,0)} f(\xi,\eta)(x-\xi)+D^{(0,1)} f(\xi,\eta)(y-\eta).\label{q1}
\fineq
Then there hold the following inequalities
$$
\Vert f-Qf\Vert_\Omega \le (1+\Vert Q\Vert)\Vert f-q^*\Vert_\Omega= (1+\Vert Q\Vert)\Vert f-q^*\Vert_T \le 
 (1+\Vert Q\Vert)\Vert f-q_1\Vert_T.
$$
By Taylor's formula, we have $f=q_1+r_1$, with
\begeq
r_1=[ D^{(1,0)}f(u,v)-D^{(1,0)} f(\xi,\eta)](x-\xi)+[\ D^{(0,1)} f(u,v)-D^{(0,1)} f(\xi,\eta)](y-\eta), \label{r1}
\fineq
the point $(u,v)$ lying somewhere in the segment joining $(\xi,\eta)$ to $(x,y)$. From that, we deduce the following upper bound
$$
\Vert f-q_1\Vert_T\le \frac{h}{2}\omega(Df, \frac{h}{2}).
$$
Actually, for the sake of simplicity, we can assume that $T$ is the triangle with vertices $(-\frac{h}{2},0), 
(0,\frac{h}{2}), (\frac{h}{2},0)$, the point $(\xi,\eta)$ being then at the origin. Due to the symmetry of the problem w.r.t. the $y$-axis, we can also assume that $(x,y)$ satisfies $x\ge 0, y\ge 0$ and $x+y\le \frac{h}{2}$.
Therefore, as the distance between $(u,v)$ and the origin is bounded above by $\frac{h}{2}$, we can write :
$$
\vert r_1(x,y)\vert \le \omega(Df, \frac{h}{2})(x+y)\le \frac{h}{2}\omega(Df, \frac{h}{2}).
$$
Finally, as $\Vert S_2\Vert\le 5$ and $\Vert W_2^*\Vert\le 3$, we obtain
$$
\Vert f-S_2f\Vert_\Omega\le 3h\, \omega(Df, \frac{h}{2})
,\;\; \Vert f-W_2^*f\Vert_\Omega\le 2h\,\omega(Df, \frac{h}{2}).
$$
\rightline{$\blacksquare$}\\


{\bf Theorem 3.} 
 {\it (i)} (Error bounds for $C^2$-functions). {\it There exists a constant $C_2>0$, with 
$C_2 \le \frac34$ for $Q=S_2$ and $C_2 \le \frac12$ for $Q=W_2^*$, such that, for $f\in C^2(\Omega)$ :}
$$
||f-Qf||_\Omega \le C_2 h^2 \om (D^2f, h/2) .
$$
 \\
 {\it (ii)} (Error bounds for $C^3$-functions). {\it There exists a constant $C_3>0$, with 
 $C_3 \le \frac18$ for $Q=S_2$ and $C_3 \le\frac{1}{12}$ for $Q=W_2^*$, such that, for  $f\in C^3(\Omega)$ :}
$$
||f-Qf||_\Omega \le C_3 h^3 ||D^3f||.
$$

{\it Proof.} 
By using a similar technique as in the proof of theorem 2, bounds on Taylor remainders can be obtained for orders 2 and 3. We have respectively, for $f\in C^2(\Omega)$ and $f\in C^3(\Omega)$  
\begeq
f=q_2+r_2 \;\; {\rm and} \;\; f=q_3+r_3,\label{q2}
\fineq
with
\begeq
r_2=\frac12 \sum_{\vert \alpha\vert=2}{2\choose \alpha}\left[ D^\alpha f(u,v)-D^\alpha f(x,y)\right]
(x-\xi)^{\alpha_1}(y-\eta)^{\alpha_2} \label{r2},
\fineq
\begeq
r_3=\frac16 \sum_{\vert \alpha\vert=3}{3\choose \alpha}D^\alpha f(u,v)
(x-\xi)^{\alpha_1}(y-\eta)^{\alpha_2}\label{r3}.
\fineq
Puting $(\xi,\eta)$ at the origin gives
$$
\vert f-q_2\vert=\vert r_2\vert\le \frac12 \omega(D^2f, \frac{h}{2})\sum_{\vert \alpha\vert=2}{2\choose \alpha}x^{\alpha_1}y^{\alpha_2}=\frac12 \omega(D^2f, \frac{h}{2})(x+y)^2\le 
\frac{h^2}{8}\omega(D^2f, \frac{h}{2}),
$$
from which we deduce respectively
$$
\Vert f-S_2f \Vert\le \frac{3h^2}{4}\omega(D^2f, \frac{h}{2})\quad {\rm and}\quad \Vert f-W_2^*f \Vert\le \frac{h^2}{2}\omega(D^2f, \frac{h}{2}).
$$
For $f\in C^3(\Omega)$, we have :
$$
\vert f-q_3\vert=\vert r_3\vert\le \frac16  \Vert D^3f\Vert \sum_{\vert \alpha\vert=3}{3\choose \alpha}x^{\alpha_1}y^{\alpha_2}=\frac16 \Vert D^3f\Vert (x+y)^3\le 
\frac{h^3}{48}\Vert D^3f \Vert,
$$
and finally we obtain
$$
\Vert f-S_2f \Vert\le \frac{h^3}{8}\Vert D^3f \Vert \quad {\rm and}\quad\Vert f-W_2^*f \Vert\le \frac{h^3}{12}\Vert D^3f \Vert.
$$
(Notice that the inequalities obtained for $W_2^*$ are better than those given in Chui-He \cite{ch} for $W_2$).
$\blacksquare$\\

{\bf Remark.} The constants of the error bounds obtained for $S_2$ are greater than the corresponding ones for $W_2^*$. These results do not mean that $S_2$ is worse than $W_2^*$, because they are a consequence of the fact that $S_2$ and $W_2^*$ belong to the same spline space, i.e. they are defined on the same triangulation  ${\cal T}_{mn}$. Therefore $N_S<N_W^*$, in particular from (\ref{quattro0}) and (\ref{sei0}) we have that $N_S=O(mn)$ and $N_W^*=O(2mn)$.

Now if we assume that  $S_2$ is defined on the triangulation  ${\cal T}_{mn}$ and that $W_2^*$ is defined on another triangulation obtained by a decomposition of $\Omega$ into ${\lceil{m\over{\sqrt{2}}}\rceil}\cdot{\lceil{n\over{\sqrt{2}}}\rceil}$ subrectangles, then the numbers of data values requested by both QIs are almost equal and the constants appearing in the respective error bounds are also comparable.

\section{Error analysis for partial derivatives}

In this section,  we compute error bounds for the first partial derivatives (Subection 5.2) of the quasi-interpolant $Qf$ in $\Omega$, and of its second partial derivatives (Subection 5.3) in the interior of each triangular cell $T_{\ell}$ of ${\cal T}_{mn}$.

\subsection{Technical lemmas}

\noi {\bf Lemma 1.} {\it Let $T$ be a triangular cell of ${\cal T}_{mn}$ included in the rectangular cell 
$\Omega_{rs}$ centered at $M_{rs}$, then: 
\begeq
\sum_{(i,j)\in K(T)}| D^\alpha B_{ij} (x,y)| \le \left\{ \begin{array}{l}
4 (h_r)^{-\alpha_1} (k_s)^{-\alpha_2} , \;{\rm for}\;\vert \alpha \vert=1 \; {\rm and}\;(x,y) \in T \\ \\ 
6 (h_r)^{-\alpha_1} (k_s)^{-\alpha_2},  \;{\rm for}\;\vert \alpha \vert=2  \; {\rm and}\; (x,y) \in int (T) \end{array} \right. \label{ventisett}
\fineq

\noi{\it Proof.}}  In the case $\vert \alpha \vert=1$, since $D^{\alpha}B_{ij}$ is a linear polynomial in the triangle $T=ABC$, we have: 
$$
|D^{\alpha}B_{ij} (x,y)| \le \max \{ |D^{\alpha}B_{ij}(A)| , \> 
|D^{\alpha}B_{ij}(B)| , |D^{\alpha}B_{ij} (C)| \}.
$$
For $\vert\alpha\vert=2$, then $D^{\alpha}B_{ij}$ is a constant inside $T$.

In \cite {s3},  the values of the first partial derivatives of $B_{ij}$ at the vertices of ${\cal T}_{mn}$ and the values of the second partial derivatives of $B_{ij}$ inside  each triangle $T$ of their support, have been computed.
Using those values, we can easily deduce the inequalities (\ref{ventisett}). \mbox{\rule{2mm}{3mm}} 

\vspace{15pt}

\noi{\bf Lemma 2.} {\it Let $Q$ be the spline operator given by (\ref{dieci}). Let $r_1$ and $r_2$ be the expressions defined by (\ref{r1}) and (\ref{r2}) for $f\in C^1(\Omega)$ and $f\in C^2(\Omega)$, respectively. 
Then for every triangle $T$ of ${\cal T}_{mn}$, the following majorations hold :
$$
\max_{(i,j)\in K(T)} |\lambda_{ij} (r_1)| \le C_1' h \om (Df, h/2),
$$
where $C_1' \le 30$ for $Q=S_2$ and $C_1'\le 35/2$ for $Q=W_2^*$;
$$
\max_{(i,j)\in K(T)} |\lambda_{ij} (r_2)| \le C_2' h^2\om (D^2 f, h/2),
$$

where $C_2'\le 61/2$ for $Q=S_2$ and $C_2'\le 65/4$ for $Q=W_2^*$.

\vspace{15pt}

\noi{\it Proof.}} We prove the desired results in the case of the  triangle $T=T_3$ (fig. 1) and we use the notations of the proof of theorem 1. For the other three types of triangles of ${\cal T}_{mn}$ we obtain the same results using a similar proof, therefore we don't report here the corresponding computations. If $Q=S_2$, then we know that
$$\begin{array}{ll}
|\mu_{ij} (r_1)| = | b_{ij} r_1(M_{ij}) + a_ir_1(M_{i-1,j}) + 
 c_i r_1(M_{i+1,j}) + \bar a_j r_1(M_{i,j-1}) + \bar c_j r_1(M_{i,j+1})| \\ \\
\le 3 |r_1(M_{i,j})|+ \frac 12 \Big[ |r_1(M_{i-1,j})| + |r_1(M_{i+1,j})| 
+ | r_1(M_{i,j-1})| + |r_1(M_{i,j+1})| \Big].
\end{array} 
$$
Moreover,  taking the origin at the midpoint $(\xi_0,\eta_0)$ of the lower edge of $T$, we can write
$$
r_1(M_{i,j})=\big[ D^{(1,0)} f(\tilde M_{i,j})- D^{(1,0)} f (O)\big] (s_i-\xi_0)
+ \big[ D^{(0,1)} f(\tilde M_{i,j}) - D^{(0,1)}f(O) \big]( t_j-\eta_0),
$$
where $\tilde M_{i,j}$ is some point lying in the segment $OM_{i,j}$.
\medskip 

Finally from the first column of Table 1, we can write
$$
\max_{(i,j)\in K(T)} |\mu_{ij}(r_1)| = \max \Big\{ | \mu_{r,s\pm1}(r_1)|, | \mu_{r\pm1,s}(r_1)|,
|\mu_{r,s}(r_1)|, |\mu_{r\pm1,s-1}(r_1)|  \Big\}=\max_{1\le k\le 7} |\mu_k(r_1)|.
$$
Recall that each central point $M_k$ has four neighbours denoted respectively $N_k,S_k,E_k,W_k$.
Here the biggest constant is obtained for $k=1$ corresponding to the central point $M_{r,s+1}$.
In that case, we have
$$
\vert \mu_1(r_1)\le 3\vert r_1(M_1)\vert+\frac12 (\vert r_1(N_1)\vert+\vert r_1(S_1)\vert+\vert r_1(E_1)\vert+\vert r_1(W_1)\vert),
$$
where, as $\Vert O\tilde M_k \Vert \le \Vert OM_k \Vert $, we can write for example
$$
\vert r_1(M_1)\vert\le \frac{3h}{2}\vert D^{(0,1)}f(\tilde M_1)-D^{(0,1)}f(O)\vert\le  \frac{3h}{2}
\omega(Df, 3\frac{h}{2})\le \frac{9h}{2} \omega(Df, \frac{h}{2}).
$$
$$
\vert r_1(W_1)\vert\le h\vert D^{(1,0)}f(\tilde W_1)-D^{(1,0)}f(O)\vert+\frac{3h}{2}\vert D^{(0,1)}f(\tilde W_1)-D^{(0,1)}f(O)\vert
$$
$$
\vert r_1(W_1)\vert\le  \frac{5h}{2}\omega(Df, \sqrt{13} \frac{h}{2})\le10\, h\,\omega(Df, \frac{h}{2}).
$$
From these inequalities and similar ones associated with the three other neighbours of $M_1$, we obtain :
$$
\vert \mu_1(r_1)\vert \le \left(\frac{27}{2}+\frac12\left(\frac{25}{2}+20+\frac12\right)\right) h\,\omega(Df, \frac{h}{2})=30\, h\,\omega(Df, \frac{h}{2}).
$$
In a similar way, we can obtain
$$
\mu_2(r_1)\;\; {\rm and}\;\; \mu_4(r_1)\le \frac{109}{4}\, h\,\omega(Df, \frac{h}{2}), \quad
\mu_3(r_1)\le \frac{17}{2} \, h\,\omega(Df, \frac{h}{2}),
$$
$$
\mu_5(r_1) \;\; {\rm and}\;\; \mu_7(r_1)\le \frac{109}{4}\, h\,\omega(Df, \frac{h}{2}), \quad
\mu_6(r_1)\le \, h\,\omega(Df, \frac{h}{2}). \;\;
$$
Finally, we obtain
$$
\max_{1\le k\le 7} |\mu_k(r_1)|\le 30\, h\,\omega(Df, \frac{h}{2}).
$$

If $Q=W^*_2$, then  from (\ref{quattro}) and (\ref{cinque}) we get
$$
| \mu_{ij}^* (r_1)| \le 2 | r_1 (M^*_{ij}) | +{1\over4}\big[ |r_1(A^*_{i-1.j-1})| + \\ \\
 |r_1(A_{i-1,j}^*)| + |r_1 (A_{i,j-1}^*)| + |r_1(A_{i,j}^*)| \big].
$$
By a procedure similar to that adopted for $S_2$, we can obtain that
$$
\max\limits_{(i,j)\in K(T)}| \mu_{ij}^* (r_1)| \le {35\over2} h\om (Df, h/2). 
$$ 
Now we consider $\lam_{ij}(r_2)$.  If $Q=S_2$, we know that
$$
|\mu_{ij} (r_2)|  \le 3 |r_2(M_{ij})| +\frac12 \big( |r_2(M_{i-1,j})| + 
 |r_2 (M_{i+1,j})| + |r_2(M_{i,j-1})| + |r_2(M_{i,j+1})| \big),
$$
where
$$
r_2(M_{i,j})=\frac12 \big\{ [D^{(2,0)} f(\tilde M_{i,j} ) -D^{(2,0)} f (O) ] (s_i-\xi_0)^2 +
 \big[ D^{(0,2)} f(M_{i,j}) -D^{(0,2)} fO)\big] (t_j-\eta_0)^2
$$
$$
 + 2 \big[ D^{(1,1)} f(M_{i,j}) - D^{(1,1)} f(O)\big] (s_i-\xi_0) (t_j-\eta_0)\big\}, 
$$
with $\tilde M_{ij}\in OM_{ij}$.
Then using a scheme similar to that  proposed for $r_1$, we get
$$
\max\limits_{(i,j)\in K(T)} |\mu_{ij} (r_2)| \le {61\over2} h^2 \om (D^2f, h) .
$$
Similarly if $Q=W_2^*$ we can deduce:
$$
\begin{array}{ll}
\max\limits_{(i,j)\in K(T)}  \mu_{ij}^* (r_2)| &\le \max\limits_{(i,j)\in K(T_1)} \big\{ 2 |r_2 (M_{ij}^*)| + 
  \frac 14 \big[ |r_2 (A_{i-1,j-1}^*)| + |r_2(A_{i-1,j}^*)|+  \\ \\ 
  & |r_2 (A_{i,j-1}^*)| + |r_2 (A_{i,j}^*)|\big]\big\} \le {65\over 4} h^2 \om (D^2 f, h/2). \qquad \quad \mbox{\rule{2mm}{3mm}}
\end{array}  
$$


\noi{\bf Lemma 3.} {\it Let $Q$ be a spline operator defined by (\ref{dieci}) and let $r_3$ be defined by (\ref{r3}) for $f\in C^3(\Om)$.  Then for every triangle $T$ of ${\cal T}_{mn}$, one has the following majoration :
\begeq
\max\limits_{(i,j)\in K(T)} | \lam_{ij} (r_3)| \le C_2'' h^3|| D^3 f || \label{quarantadu} \fineq
with $C_2'' \le 269/48$ for $Q=S_2$, $C_2'' \le 65/24$ for $Q=W_{2}^*$.}

\vspace{15pt}

\noi {\it Proof.} As in the proof of Lemma 2, we only discuss here the case of a triangle $T=T_3$ and we don't report the computations for the  three other types of triangles.

If $Q=S_2$ then, from (\ref{r3}), we can write
\begeq
r_3(M_{i,j}) = \frac16 \sum_{\vert \alpha\vert=3}{3\choose \alpha}D^\alpha f({\tilde M_{i,j}})
(s_i-\xi_0)^{\alpha_1}(t_j-\eta_0)^{\alpha_2}, \label{quarantatre} 
\fineq
\\
with $\tilde M_{i,j}\in OM_{i,j}$.
Therefore if we denote
\begeq
\varphi(s_i,t_j) =|\sum_{\vert \alpha\vert=3}{3\choose \alpha}
(s_i-\xi_0)^{\alpha_1}(t_j-\eta_0)^{\alpha_2}| \label{quarantaquattr}
\fineq
\\
and if we proceed as in the proof of  Lemma 2, from (\ref{due}), (\ref{quarantatre}), (\ref{quarantaquattr}) we can deduce the desired result. Indeed:
$$
\begin{array}{ll}
 \max\limits_{(i,j)\in K(T)} |\mu_{ij}(r_3)| \le{1\over6}|| D^3 f|| \big\{ 3\varphi (s_r, t_{s+1}) + {1\over2} [ \var (s_{r-1}, t_{s}) + \var (s_{r+1},t_{s+1}) \\ \\
+ \var (s_r, t_{s+2})+ \var (s_r, t_s) \big] \big\} \le 
{1\over6} \Del^3 || D^3 f || \big\{ {1\over2} [ 3({5\over2})^3 + {1\over8}] +3({5\over2})^3 \big\} 
 \le {269\over48} \Del^3 ||D^3 f||. \end{array} 
$$
Using the same method, we can prove (\ref{quarantadu}) for $Q=W_2^*$.{\mbox{\rule{2mm}{3mm}}

\subsection{Error estimates on first partial derivatives}


\noi{\bf Theorem 4.}  ($C^1$ functions). {\it For $\vert \alpha \vert=1$,  there exists a constant $\ov C_1>0$, with $\ov C_1\le 120$ for $Q=S_2$ and $\ov C_1 \le 70$ for $Q=W_2^*$, such that, for $f\in C^1(\Om)$}

\begeq
||D^{\alpha} f - D^{\alpha} Qf||_\Omega \le \Big[1+\ov C_1 \Big({h\over\del}\Big)\Big] \om (Df, h/2) \label{quarantasei}
\fineq

\noi{\it Proof.} For $\vert\alpha\vert=1$ we consider a closed triangular cell $T$ of ${\cal T}_{mn}$ where
$$
||D^{\alpha} f - D^{\alpha} Qf||_\Omega = ||D^{\alpha} f- D^{\alpha} Qf||_{T}. 
$$

For any point $P=(x,y)\in T$, since $Q$ reproduces $\PP_2$, we can write: 
\begeq
|D^{\alpha} f(P) - D^{\alpha} Qf(P)| \le |D^{\alpha} f(P)  - D^{\alpha} q_1(P) | + 
|D^{\alpha} Q(f-q_1)(P) |, \label{quarantotto}
\fineq
\\
with $q_1$ defined in (\ref{q1}).

We remark that, from (\ref{q1}), there results:
\begeq
|D^{\alpha} f(P)-D^{\alpha} q_1(P) | = | D^{\alpha} f(P) - D^{\alpha} f(O)| \le \om (Df, h/2)\label{dalfa}
\fineq
and
$$
|D^{\alpha} Q(f-q_1)(P)| \le \max\limits_{(i,j)\in K(T)} |\lam_{ij}(r_1)|  \sum\limits_{(i,j)\in K(T)} |D^{\alpha} B_{ij}(P)| 
$$
with $r_1$ defined in (\ref{r1}). \medskip 

Moreover we recall that, from Lemma 1, for $\vert\alpha\vert=1$
\begeq
\sum\limits_{(i,j)\in K(T)} |D^{\alpha} B_{ij}(P)| \le 4 \del^{-2}. \label{cinquantuno}
\fineq
\medskip 

Finally from Lemma 2 and (\ref{cinquantuno}) we obtain 
\begeq
|D^{\alpha} Q(f-q_1)| \le 4 C_1' \, {h\over\del} \om (Df, h/2) \label{cinquantadue}
\fineq

Therefore, from (\ref{quarantotto}), (\ref{dalfa}) and (\ref{cinquantadue}), the result (\ref{quarantasei}) follows with $\ov C_1=4C_1'$. \rule{2mm}{3mm}.

\vspace{15pt}

\noi{\bf Theorem 5.} (i)  ($C^2$ functions). {\it There exists a constant $\bar C_2>0$, with  $\bar C_2\le 122$  for $Q=S_2$ and $\bar C_2\le 65$ for $Q=W_2^*$, such that, for $f\in C^2(\Omega)$ and $\vert\alpha\vert=1$ :}
\begeq
|| D^{\alpha} f- D^{\alpha}  Qf||_\Omega \le \left[ 1 +  \bar C_2 \left( \frac h \del \right) \right] h \om (D^2f, h/2) \label{cinquantatre}
\fineq

(ii)  ($C^3$ functions).  {\it  Then there exists a constant $\bar C_3>0$, with $\bar C_3\le 
{269\over12}$ for $Q=S_2$ and $\bar C_3\le {65\over 6}$  for $Q=W_2^*$, such that,  for  $f\in C^3(\Omega)$ and $\vert\alpha\vert=1$ :
\begeq
||D^{\alpha}f - D^{\alpha} Qf|| \le \bar C_3
\left( \frac h \del \right) h^2|| D^3f|| \label{cinquantaquattro}
\fineq}

\noi{\it Proof.} The proof is similar to that of Theorem 4. \vspace{15pt}

For $\vert\alpha\vert=1$ and  $P \in T$, we can write:
\begeq
|D^{\alpha}f(P) - D^{\alpha} Qf(P)| \le | D^{\alpha}f(P) - D^{\alpha} q_2(P) | + | D^{\alpha} Q(f-q_2)(P) | \label{cinquantacinque}
\fineq
where $q_2 \in \PP_2$ has been defined in (\ref{q2}).

We remark that
\begeq
| D^{\alpha}f(P) - D^{\alpha} q_2(P) | \le  h\om (D^2f, h/2). \label{cinquantasei}
\fineq

Moreover
\begeq
| D^{\alpha} Q(f-q_2)(P) | \le \max\limits_{(i,j)\in K(T)}|\lam_{ij}(r_2)|  \sum_{(i,j)\in K(T)} | D^{\alpha} B_{ij}(P)| \label{cinquantasette},
\fineq
with $r_2$ defined in (\ref{r2}).

Now from (\ref{cinquantasette}), Lemma 1 and Lemma 2 we can write that
\begeq
|D^{\alpha} Q(f-q_2)(P)| \le 4 C_2' \left( \frac h \del \right) h \om (D^2f, h/2). \label{cinquantotto}
\fineq
Finally from (\ref{cinquantacinque}), (\ref{cinquantasei}) and (\ref{cinquantotto}), the result (\ref{cinquantatre}) follows, with ${\ov C}_2=4C_2'$. \medskip\medskip

Moreover if $f\in C^3(\Omega)$, from (\ref{r3}), there results: 
\begeq 
|D^{\alpha}(f-q_3)(P)| \le \frac {h^2}2 || D^3 f|| \label{cinquantanove},
\fineq
and from Lemma 1 and Lemma 3
\begeq
|D^{\alpha} Q(f-q_3)(P)\vert \le 4 C_2'' \left( \frac h \del \right) h^2 ||D^3f ||. \label{sessanta}
\fineq

Therefore from (\ref{cinquantacinque}), (\ref{cinquantanove})  and (\ref{sessanta}) the result (\ref{cinquantaquattro}) follows, with $\bar C_3 = 4C_2''$. \rule{2mm}{3mm}


\subsection{Error estimates on second partial derivatives}


\noi{\bf Theorem 6.} 
(i) ($C^2$ functions). {\it There exists a constant  $D_2>0$, with $D_2\le {183}$ for $Q=S_2$, and $D_2\le {195\over2}$ for $Q=W_2^*$,  such that for $f\in C^2(\Omega)$ and  $\vert \alpha\vert=2$ :}
\begeq
||D^{\alpha}f - D^{\alpha} Qf ||_{{\rm int}{(T)}} \le \Big[ 1+ D_2 
\left( \frac h \del \right)^2 \Big] \om (D^2f, h/2) \label{sessantuno}.
\fineq

(ii) ($C^3$ functions). { \it There exists a constant  $D_3>0$, with $D_3 \le {269\over8}$ if $Q=S_2$ and $D_3 \le {65\over4}$ if $Q=W_2^*$, such that, for  $f\in C^3(\Omega)$ :}
\begeq
||D^{\alpha}f - D^{\alpha} Qf ||_{{\rm int}(T)} \le \Big[ 1 + D_3 \left( \frac h \del \right)^2 \Big] h ||D^3 f|| \label{sessantadue}.
\fineq}

\noi{\it Proof.} For $\vert \alpha\vert=2$  and for any $P\in$ int$(T)$,  we have 
\begeq
|D^{\alpha}f - D^{\alpha} Qf | \le | D^{\alpha}f - D^{\alpha} q_2| + | D^{\alpha} Q(f-q_2)|, \label {sessantatre}
\fineq
with $q_2$ defined in (\ref{q2}).

From (\ref{q2}), we deduce 
\begeq
| D^{\alpha}f - D^{\alpha} q_2| \le \om (D^2f, h/2). \label{sessantaquattro}
\fineq

Moreover
\begeq
| D^{\alpha} Q(f-q_2)| \le \max\limits_{(i,j)\in K(T)} |\lam_{ij}(r_2)| \cdot \sum_{(i,j)\in K(T)} | D^{\alpha} B_{ij}|, \label{sessantacinque}
\fineq
with $r_2$ defined in (\ref{r2}).

From (\ref{sessantacinque}), Lemma 1 and Lemma 2 we obtain that  for $P\in$int$(T)$:
\begeq
| D^{\alpha} Q(f-q_2)(P)| \le 6  C_2' \left( \frac h \del \right)^2 \om (D^2f, h/2) \label {sessantasei}
\fineq
Therefore from (\ref{sessantatre}), (\ref{sessantaquattro}) and (\ref{sessantasei}) the result (\ref{sessantuno}) follows, with $D_2 = 6C_2'$.

Finally if $f\in C^3(\Omega)$, for $P\in$ int$(T)$ and $\vert \alpha\vert=2$, from (\ref{sessantatre}), Lemma 1 and Lemma 3, we obtain the result (\ref{sessantadue}), with $D_3 = 6 C_2''$. \rule {2mm}{3mm}\\

We note that the same remarks given at the end of Section 3 are also valid  for the error bound constants of the above theorems.

 
\subsection{Convergence for quasi-uniform partitions}

Assume that the sequence of partitions $\{X_m \times Y_n\}$ of $\Omega$ is $\gam$-{\sl quasi uniform} i.e. there exists a  constant $\gam>1$ such that 
$$
0<  h_{mn}/\delta_{mn}\le \gamma,
$$
where $h_{mn}$ and $\delta_{mn}$ are respectively the maximum and the minimum steplengths of the partition $\{X_m \times Y_n\}$. Then the following theorem shows that for both $Q=S_2$ and $Q=W_2^*$ 
$$
D^{\alpha} Qf \rightarrow D^{\alpha} f \qquad {\rm as} \quad  h_{mn} \to 0 
$$
in $\Omega$ for $\vert\alpha\vert=1$, and in the interior of each triangular cell $T$ of ${\cal T}_{mn}$, for $\vert\alpha\vert=2$.\\

\noi{\bf Theorem 7.} {\it Let $\{ X_m \times Y_n\}$ a $\gam$-quasi uniform sequence of partitions.  

{\em (i)} If $f\in C^s(\Omega)$, $s=1,2$, then for $\vert\alpha\vert=1$
$$
|| D^{\alpha}f - D^{\alpha} Qf ||_\Omega = O (h_{mn}^{s-1} \om (D^sf, h_{mn}/2))
$$
{\em (ii)} If $f \in C^2(\Omega)$  then for}$\vert\alpha\vert=2$
$$
|| D^{\alpha}f - D^{\alpha} Qf ||_{{\rm int}(T)} = O (\om (D^2 f, h_{mn}/2)) 
$$
\vspace{10pt}
\noi {\it Proof.} The result immediately follows from the $\gam$-quasi \ uniformity \ of \ $\{X_m \times Y_n\}$ and from Theorems 4, 5, 6. \rule {2mm}{3mm} \\


\end{document}